# Boundary integrated neural networks (BINNs) for acoustic radiation and scattering

Wenzhen Qu[1], Yan Gu[1*], Shengdong Zhao[1], Fajie wang[2]

[1]*School of Mathematics and Statistics, Qingdao University, Qingdao 266071, PR China*
[2]*College of Mechanical and Electrical Engineering, National Engineering Research Center for Intelligent Electrical Vehicle Power System, Qingdao University, Qingdao 266071, China*

**Abstract**

This paper presents a novel approach called the boundary integrated neural networks (BINNs) for analyzing acoustic radiation and scattering. The method introduces fundamental solutions of the time-harmonic wave equation to encode the boundary integral equations (BIEs) within the neural networks, replacing the conventional use of the governing equation in physics-informed neural networks (PINNs). This approach offers several advantages. Firstly, the input data for the neural networks in the BINNs only require the coordinates of "boundary" collocation points, making it highly suitable for analyzing acoustic fields in unbounded domains. Secondly, the loss function of the BINNs is not a composite form, and has a fast convergence. Thirdly, the BINNs achieve comparable precision to the PINNs using fewer collocation points and hidden layers/neurons. Finally, the semi-analytic characteristic of the BIEs contributes to the higher precision of the BINNs. Numerical examples are presented to demonstrate the performance of the proposed method.

---

[*] Corresponding author, Email: guyan1913@163.com (Y. Gu)

# 1. Introduction

The boundary element method (BEM) has gained recognition as a formidable technique for numerically analyzing acoustic fields, owing to its semi-analytical nature and boundary-only discretization [1, 2]. By incorporating fundamental solutions into the BEM, the time-harmonic wave equation for acoustic problems, along with boundary conditions and the Sommerfeld radiation condition at infinity, can be transformed into boundary integral equations (BIEs) [3]. Consequently, the BEM offers several advantages, including the reduction of problem dimensionality by one and the direct solution of unbounded domain problems without the need for special treatments.

Over the past decade, significant attention has been directed towards machine learning, owing to the remarkable advancements in computing resources and the abundance of available data [4]. Among the prominent tools in machine learning, deep neural networks (DNNs) have emerged as outstanding approximations of functions, demonstrating immense potential for numerical simulations of partial differential equation (PDE) problems [5, 6]. Up to now, numerous DNN-based approaches have been devised to tackle PDEs, including physics-informed neural networks (PINNs) [7-9], the deep Galerkin method (DGM) [10, 11], and the deep Ritz method (DRM) [12, 13]. The aforementioned DNN-based methods directly approximate the solution of problems using a neural network. Subsequently, a loss function or composite form is constructed, incorporating information from the residuals of the partial differential equation (PDE) with boundary/initial conditions or the energy functional form.

There have been remarkable contributions in acoustic numerical analysis through the utilization of DNN-based methods [14, 15]. The DNNs are typically trained and applied in finite domains, which poses challenges when directly using them to solve unbounded domain problems. Very recently, Lin et al. [16] made the first attempt to integrate neural networks with indirect boundary integral equations (BIEs) for solving partial differential equation (PDE) problems with Dirichlet boundary conditions. Following this, Zhang et al. [17] utilized neural networks to approximate solutions of direct BIEs using non-uniform rational B-splines (NURBS) parameterization of the boundary for potential problems. The aforementioned approaches are theoretically well-suited for addressing problems in unbounded domains. However, they have not been empirically validated by related problems in the references mentioned.

In this paper, we propose a novel approach called the boundary integrated neural networks (BINNs) to analyze acoustic problems in both bounded and unbounded domains. The method involves the approximated solutions of neural networks, trained solely on boundary collocation points, into the direct acoustic boundary integral equations (BIEs) using quadratic elements. The loss function is then constructed based on the BIE residuals and is minimized specifically at these collocation points. Three numerical examples with various types of boundary conditions are provided to validate the proposed method. The numerical results obtained using the developed approach are compared with those obtained using the PINNs as well as the exact solutions.

## 2. Mathematical formulation for acoustic problem

The time-harmonic wave equation [18], commonly referred to as the Helmholtz equation, can be expressed in 2D domain $\Omega$ as follows:

$$\nabla^2 p + k^2 p = 0, \ p \in \Omega \tag{1}$$

where $p$ represents the complex acoustic pressure, while $k$ denotes the wave number. The wave number, defined as $\omega/c$, corresponds to the ratio of the angular frequency $\omega$ to the speed of the acoustic wave $c$ in the medium $\Omega$. The equation (1) is subject to Dirichlet and Neumann boundary conditions (BCs) as

$$p(\boldsymbol{x}) = \overline{p}(\boldsymbol{x}), \ \boldsymbol{x} \in \Gamma_D, \tag{2}$$

$$q(\boldsymbol{x}) = \frac{\partial \overline{p}(\boldsymbol{x})}{\partial \boldsymbol{n}(\boldsymbol{x})} = \mathrm{i}\rho\omega\overline{v}(\boldsymbol{x}), \ \boldsymbol{x} \in \Gamma_N, \tag{3}$$

where $\boldsymbol{n}(\boldsymbol{x})$ represents the outward unit normal vector to the boundary $\Gamma$ at point $\boldsymbol{x}$, $\rho$

denotes the density of the medium, i means the imaginary unit, $v(x)$ is the normal velocity, and the upper bars on the pressure and normal velocity indicate the known functions. Furthermore, as the distance $r$ from the source tends to infinity, it is essential for the pressure field to satisfy the Sommerfeld radiation condition as

$$\lim_{r \to \infty} \sqrt{r}\left(\frac{\partial p(r)}{\partial r} - \mathrm{i}kp(r)\right) = 0 \qquad (4)$$

## 3. Boundary integrated neural networks (BINNs)

*3.1. Boundary integral equations (BIEs)*

By incorporating the fundamental solutions, the time-harmonic wave equation for acoustic pressure can be transformed into an integral form [3], represented as

$$p(\boldsymbol{x}) + \int_{\Gamma} F(\boldsymbol{x}, \boldsymbol{y}) p(\boldsymbol{y}) d\Gamma(\boldsymbol{y}) = \int_{\Gamma} G(\boldsymbol{x}, \boldsymbol{y}) q(\boldsymbol{y}) d\Gamma(\boldsymbol{y}), \quad \boldsymbol{x} \in \Omega \qquad (5)$$

where $\boldsymbol{x}$ and $\boldsymbol{y}$ represent the source and field points, respectively, while $G(\boldsymbol{x}, \boldsymbol{y})$ and $F(\boldsymbol{x}, \boldsymbol{y})$ respectively denote the fundamental solutions of the time-harmonic wave equation and its corresponding normal derivatives. $G(\boldsymbol{x}, \boldsymbol{y})$ and $F(\boldsymbol{x}, \boldsymbol{y})$ for 2D problems are defined as

$$G(\boldsymbol{x}, \boldsymbol{y}) = \frac{\mathrm{i}}{4} H_0^{(1)}\left(kr(\boldsymbol{x}, \boldsymbol{y})\right), \text{ and } F(\boldsymbol{x}, \boldsymbol{y}) = \frac{\partial G(\boldsymbol{x}, \boldsymbol{y})}{\partial n(\boldsymbol{y})} \qquad (6)$$

where $H_0^{(1)}$ represents the first kind Hankel function of order zero, $r$ is the distance between points $\boldsymbol{x}$ and $\boldsymbol{y}$. Taking the limit as $\boldsymbol{x}$ in Eq. (5) approaches the boundary $\Gamma$, we obtain

$$C(\boldsymbol{x})p(\boldsymbol{x}) + \int_{\Gamma}^{\mathrm{CPV}} F(\boldsymbol{x}, \boldsymbol{y}) p(\boldsymbol{y}) d\Gamma(\boldsymbol{y}) = \int_{\Gamma} G(\boldsymbol{x}, \boldsymbol{y}) q(\boldsymbol{y}) d\Gamma(\boldsymbol{y}), \quad \boldsymbol{x} \in \Omega \qquad (7)$$

in which $C(\boldsymbol{x}) = 0.5$ as the boundary near point $\boldsymbol{x}$ is smooth, and $\int_{\Gamma}^{\mathrm{CPV}}$ denotes the integral evaluated in the sense of Cauchy principal value (CPV). In this study, regular integrals are computed using the standard Gaussian quadrature with twenty Gaussian points, while the singular integrals are evaluated using a direct method developed by Guiggiani and Casalini [19] for CPV integrals. It is widely acknowledged that the handling techniques for singular integrals in boundary integral equations (BIEs) have reached a high level of maturity. However, the detailed methods for handling singular integrals are beyond the scope of this work. Interested readers are referred to relevant references for further information.

*3.2. Discretization of BIEs*

We discretize the BIEs using discontinuous quadratic element [20]. The shape functions, denote as $N_i(\xi)(i=1,2,3)$, of the elements are assumed to have the following forms:

$$N_1(\xi) = \frac{\xi(\xi-1)}{2}, \quad N_2(\xi) = (1-\xi)(1+\xi), \quad \text{and} \quad N_3(\xi) = \frac{\xi(\xi+1)}{2} \tag{8}$$

in which $\xi \in [-1,1]$ indicates the dimensionless coordinate. Then, the geometry of each quadratic element can be described as

$$\boldsymbol{y} = N_1(\xi)\boldsymbol{y}_1 + N_2(\xi)\boldsymbol{y}_2 + N_3(\xi)\boldsymbol{y}_3 \tag{9}$$

where $\boldsymbol{y}_1(\xi=-1), \boldsymbol{y}_2(\xi=0),$ and $\boldsymbol{y}_3(\xi=1)$ denote the right, middle, and left points of the mentioned boundary element as shown in Fig. 1, respectively. The pressure and its normal derivative on the boundary element are approximated by quantities $p_i, q_i (i=1,2,3)$ on points $\boldsymbol{y}'_1(\xi=-\alpha), \boldsymbol{y}'_2(\xi=0),$ and $\boldsymbol{y}'_3(\xi=\alpha)$ in Fig. 1, expressed as follows:

$$p(\boldsymbol{y}) = N_1\left(\frac{\xi}{\alpha}\right)p_1 + N_2\left(\frac{\xi}{\alpha}\right)p_2 + N_3\left(\frac{\xi}{\alpha}\right)p_3 \tag{10}$$

$$q(\boldsymbol{y}) = N_1\left(\frac{\xi}{\alpha}\right)q_1 + N_2\left(\frac{\xi}{\alpha}\right)q_2 + N_3\left(\frac{\xi}{\alpha}\right)q_3 \tag{11}$$

where $\alpha \in (0,1)$. In the numerical calculations of this work, the value of $\alpha$ is set to 0.8, and its influence on the numerical results is negligible.

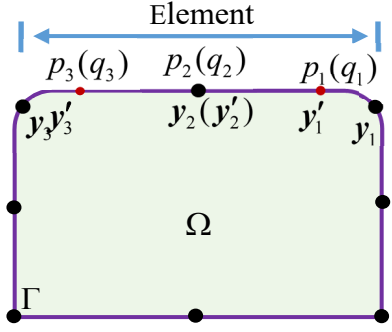

**Figure 1.** Discontinuous quadratic element

Based on the aforementioned discontinuous quadratic element, the discretized form of the BIE (7) is given as

$$C(\boldsymbol{x}^m)p(\boldsymbol{x}^m)+\sum_{i=1}^{N}\sum_{j=1}^{3}p_j^i\int_{-1}^{1}F(\boldsymbol{x}^m,\boldsymbol{y}_i(\xi))N_j(\xi/\alpha)J_i(\xi)d\xi=\sum_{i=1}^{N}\sum_{j=1}^{3}q_j^i\int_{-1}^{1}G(\boldsymbol{x}^m,\boldsymbol{y}_i(\xi))N_j(\xi/\alpha)J_i(\xi)d\xi \quad (12)$$

where $N$ represents the number of boundary elements, $\boldsymbol{x}^m (m=1,2,...,3N)$ are boundary collocation points and selected to be the same set as points $\boldsymbol{y}_j'$ (see Fig. 1) on these elements, $p_j^i$ and $q_j^i$ respectively denote the pressure and its normal derivative at the $j$-th collocation point of the $i$-th element, and $J_i(\xi)$ represents the Jacobian of transformation from the global coordinate $\boldsymbol{y}$ to the dimensionless coordinate $\xi$ for integrals at the $i$-th element.

After discretizing the BIEs through the process mentioned earlier, we can define the following two functions with Eq. (12) to facilitate the construction of the loss function in subsequent steps

$$LE(\boldsymbol{x}^m,\boldsymbol{p})=C(\boldsymbol{x}^m)p(\boldsymbol{x}^m)+\sum_{i=1}^{N}\sum_{j=1}^{3}p_j^i\int_{-1}^{1}F(\boldsymbol{x}^m,\boldsymbol{y}_i(\xi))N_j(\xi/\alpha)J_i(\xi)d\xi \quad (13)$$

$$RE(\boldsymbol{x}^m,\boldsymbol{q})=\sum_{i=1}^{N}\sum_{j=1}^{3}q_j^i\int_{-1}^{1}G(\boldsymbol{x}^m,\boldsymbol{y}_i(\xi))N_j(\xi/\alpha)J_i(\xi)d\xi \quad (14)$$

where $\boldsymbol{p}=\{p_j^i\}_{j=1,2,3}^{i=1,...,N}$ and $\boldsymbol{q}=\{q_j^i\}_{j=1,2,3}^{i=1,...,N}$

*3.3. Neural networks and loss function of the BINNs*

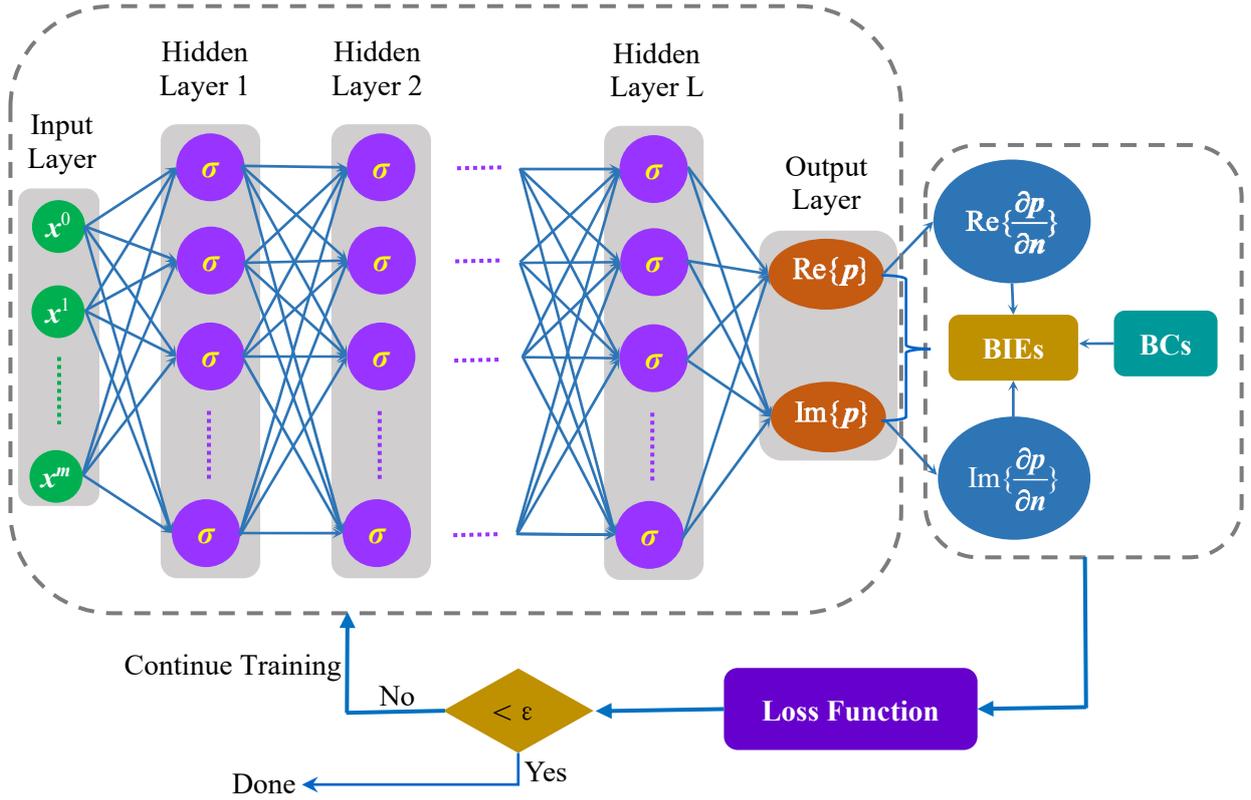

**Figure 2.** The framework of the boundary integrated neural networks (BINNs).

We present the construction of the BINNs by seamlessly integrating neural networks and the BIEs in this subsection. As illustrated in Fig. 2, we utilize a fully connected neural architecture including the input layer, the $L$ hidden layers, and the output layer. The number of neurons in $l$ hidden layer is set to $n_l$. Based on the neural networks approximation, the real and imaginary parts of trial solutions of pressures at a collocation point $x$ can be expressed as

$$\text{Re}\{p(\boldsymbol{x},\boldsymbol{w},\boldsymbol{b})\} = \hbar_1\left(\lambda_L\left(\lambda_{L-1}\left(...\left(\lambda_1(\boldsymbol{x})\right)\right)\right)\right) \tag{15}$$

$$\text{Im}\{p(\boldsymbol{x},\boldsymbol{w},\boldsymbol{b})\} = \hbar_2\left(\lambda_L\left(\lambda_{L-1}\left(...\left(\lambda_1(\boldsymbol{x})\right)\right)\right)\right) \tag{16}$$

where $\hbar_k (k=1,2)$ and $\lambda_l (l=1,2,...,L)$ are linear and nonlinear mappings, expressed as follows

$$\hbar_k(g) = \boldsymbol{w}'_k * g + b'_k \tag{17}$$

$$\lambda_l(g) = \sigma(\boldsymbol{w}_l * g + \boldsymbol{b}_l) \tag{18}$$

with weights $\boldsymbol{w}'_k \in R^{n_L}, \boldsymbol{w}'_l \in R^{n_l * n_{l-1}}$ ($n_0 = 2$), biases $b'_k \in R, \boldsymbol{b}'_l \in R^{n_l}$, and the activation function $\sigma$.

Here, Table 1 lists some commonly used activation functions. To obtain the normal derivatives of acoustic pressures approximated by the above neural networks, we employ the "dlgradient", which is an automatic differentiation function in the Deep Learning Toolbox of MATLAB.

Table 1. Some commonly used activation functions.

| | Arctan | Sigmoid | Swish | Softplus | Tanh |
|---|---|---|---|---|---|
| $\sigma(z)$ | $\arctan(z)$ | $\dfrac{1}{1+e^{-z}}$ | $\dfrac{z}{1+e^{-z}}$ | $\ln(1+e^{z})$ | $\dfrac{e^{z}-e^{-z}}{e^{z}+e^{-z}}$ |
| | 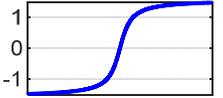 | 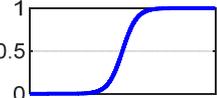 | 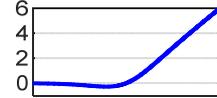 | 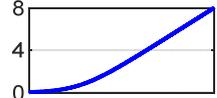 | 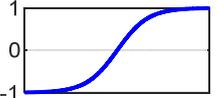 |

We construct two different forms of loss functions and will explore their performance in next sections. Firstly, we incorporate the known pressures $p$ and/or normal derivatives $q$ directly into the BIEs, creating the following loss function referred to as $Loss$

$$Loss = \frac{1}{3N}\sum_{m=1}^{3N}\left(LE(\boldsymbol{x}^m,\boldsymbol{p}) - RE(\boldsymbol{x}^m,\boldsymbol{q})\right)^2 \tag{19}$$

where the unknown $p$ and/or $q$ on the boundary are approximated by the neural networks. For the second form of the loss function, we approximate both the pressures and normal derivatives in the BIEs and boundary constraints using the neural networks. The loss function named as Loss $Loss_{BC}$ is then constructed as

$$Loss_{BC} = \frac{1}{3N}\sum_{m=1}^{3N}\left(LE(\boldsymbol{x}^m,\boldsymbol{p}) - RE(\boldsymbol{x}^m,\boldsymbol{q})\right)^2 + \frac{1}{N_D}\sum_{m=1}^{N_D}\left(\boldsymbol{p}_D - \overline{\boldsymbol{p}}_D\right)^2 + \frac{1}{N_N}\sum_{m=1}^{N_N}\left(\boldsymbol{q}_N - \overline{\boldsymbol{q}}_N\right)^2 \tag{20}$$

where the subscripts $D$ and $N$ respectively denote the Dirichlet and the Neumann BC, $N_D$ and $N_N$ indicate the numbers of Dirichlet and Neumann boundary collocation points respectively, and the superscript bar represents the known quantities.

*3.4. Optimization of parameters and solution of pressure at interior point*

In the previous subsections, we have established the architecture of the neural networks and defined the loss function for the BINNs. The next step is to optimize the weights and biases of each neuron by minimizing the corresponding loss function, either Eq. (19) or Eq. (20). To accomplish this optimization process, we utilize the powerful and widely used "fmincon" function in MATLAB. The "fmincon" is specifically designed to minimize constrained nonlinear multivariable functions.

By applying this optimization approach, we are able to obtain accurate numerical results for the unknown pressures and normal derivatives along the boundary. Once the pressures $p$ and normal derivatives $q$ at the boundary collocation points are determined, we can easily calculate the numerical solution for the pressure at any interior point using Eq. (5).

## 4. Numerical examples

To evaluate the performance of the BINNs, several benchmark examples involving bounded and unbounded domains under various BCs are provided. The accuracy of the present approach is thoroughly investigated by examining the influence of parameters such as the hidden layer number, neuron number in each layer, and the choice of activation function. The numerical results calculated by the BINNs are compared against those obtained using the traditional PINNs as well as the theoretical solutions.

All the MATLAB codes used in this study are executed on a computer equipped with an Intel Core i9-11900F 2.5 GHz CPU and 64 GB of memory. The precision of the numerical results is assessed using relative error, which is defined as

$$\text{Relative error (RE)} = \sqrt{\sum_{i=1}^{M}(\tilde{p}_i - p_i)^2} \Big/ \sqrt{\sum_{i=1}^{M} p_i^2} \qquad (21)$$

where $M$ denote the number of calculated points, $\tilde{p}_i$ and $p_i$ are numerical and analytical solutions at $i$-th calculated point, respectively.

*4.1. Interior acoustic field*

As the first example, we consider the distribution of acoustic pressure in a rectangle domain with a length of 3 m and a height of 1.5 m, as illustrated in Fig. 3. The center of the domain is $(1.5, 0.75)$. The boundary is subject to two different cases of BCs.

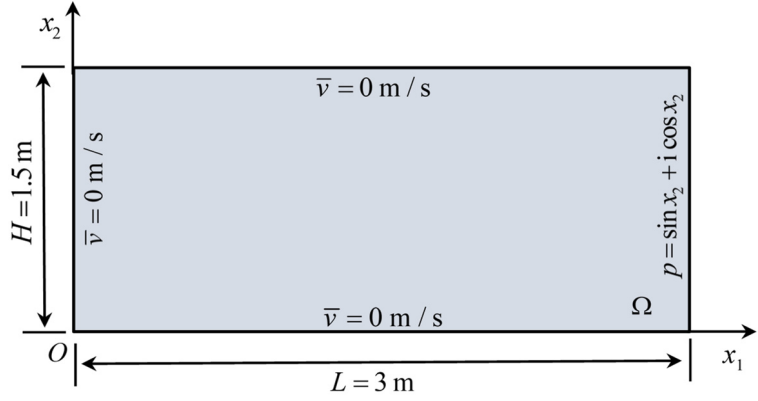

**Figure 3.** The dimension of the rectangle domain and the BCs of case 2.

Case 1: Dirichlet BC

The pressure on the boundary is specified as

$$p(x_1', x_2') = \cos(kx_1') + i \sin(kx_2'), \quad (x_1', x_2') \in \Gamma \qquad (22)$$

Obviously, the analytical solution for this case is $p(x_1, x_2) = \cos(kx_1) + i \sin(kx_2)$, $(x_1, x_2) \in \Omega$.

Initially, we assess the performance of the BINNs using two distinct forms of loss functions. Four distinct neural architectures are configured as follows: a) a single hidden layer consisting of 10 neurons; b) a single hidden layer consisting of 20 neurons; c) two hidden layers, each with 10 neurons; d) two hidden layers, each with 20 neurons. The training process for optimization stops when the iteration count reaches 10000. A total of 270 boundary collocation points, corresponding to 90 boundary elements, are utilized. The activation function selected for neural networks is $\sigma(z) = z / (1 + e^{-z})$. The wave number is set to $k = 2$ m$^{-1}$.

Using the BINNs with $Loss$ and $Loss_{BC}$, Table 2 presents the relative errors of the real and imaginary components of the pressure along the evaluated line $x_2 = 0.75\,\text{m}$ with 30 equally spaced points for calculation purposes. The numerical results obtained through the use of $Loss$ showcase superior accuracy when compared to the results obtained using $Loss_{BC}$. Remarkably, even using the networks with a single hidden layer consisting of 10 neurons, the present method with $Loss$ achieves high accuracy in the numerical results. Additionally, there is a slight improvement when employing more hidden layers or increasing the number of neurons in each layer. In contrast, the BINNs with $Loss_{BC}$ requires a greater number of hidden layers and neurons to attain sufficiently accurate results.

**Table 2.** Errors of pressures by the BINNs with $Loss$ and $Loss_{BC}$ based on four neural networks

| Error | $Loss$ | | | | $Loss_{BC}$ | | | |
|---|---|---|---|---|---|---|---|---|
| | a | b | c | d | a | b | c | d |
| $\text{Re}\{p\}$ | 2.38E-06 | 4.66E-07 | 5.74E-07 | 1.38E-07 | 2.07E-02 | 2.77E-03 | 4.09E-05 | 3.85E-05 |
| $\text{Im}\{p\}$ | 6.43E-07 | 7.30E-08 | 2.88E-07 | 9.12E-08 | 3.61E-03 | 7.41E-04 | 8.96E-05 | 5.62E-05 |

Fig. 4 illustrates the convergence process of two designated loss functions $Loss$ and $Loss_{BC}$ over iterations ranging from 1 to 10000, with values recorded at every 100 iterations. It is apparent that $Loss$ exhibits a faster convergence rate compared to $Loss_{BC}$. Therefore, the BINNs with $Loss$ has a better performance in comparison to that with $Loss_{BC}$, as indicated in Table 2. To expedite the convergence process of the loss function $Loss_{BC}$, the incorporation of additional learning techniques is necessary to balance its different loss terms. Consequently, $Loss$ stands as the superior choice for an efficient loss function in the context of BINNs when compared to $Loss_{BC}$. Henceforth, the BINNs will employ $Loss$ in all subsequent computational processes unless otherwise specified.

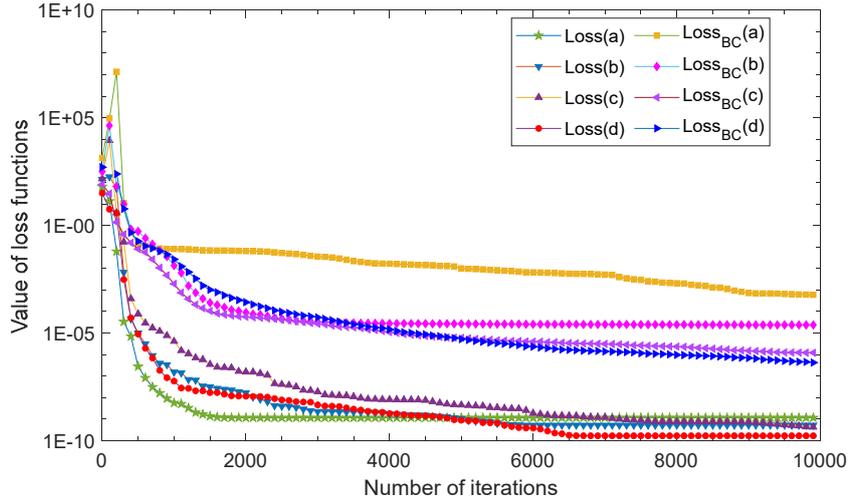

**Figure 4.** Convergence process of loss functions constructed with different neural architectures

Next, we present a comparison of the accuracy of the numerical results obtained using the BINNs and the traditional PINNs. The same calculated points are distributed on the line $x_2 = 0.75 \text{ m}$. The wave number, activation functions of the neural networks, and optimization stopping criteria for both methods remain consistent with the previous settings. The BINNs adopts a single hidden layer comprising 20 neurons, while the PINNs utilizes two different networks: a) a single hidden layer with 20 neurons, and b) three hidden layers, each containing 20 neurons. The collocation points for the PINNs are uniformly distributed within the rectangular domain and its boundary, while for the BINNs, they are only placed on the boundary. Fig. 5 and Fig. 6 plot the convergence curves of the pressures obtained by the BINNs and the PINNs as the number of collocation points increases. Clearly, the BINNs exhibits a faster and more stable convergence rate compared to the PINNs with networks "a" or "b". Furthermore, the precision of the pressures evaluated by the BINNs is higher even with a smaller number of collocation points, as compared to the PINNs. Therefore, in order to achieve comparable precision in pressure calculations, the BINNs necessitates significantly fewer collocation points and hidden layers/neurons compared to the PINNs. This observation also demonstrates that the BINNs exhibits higher computational efficiency in comparison to the PINNs.

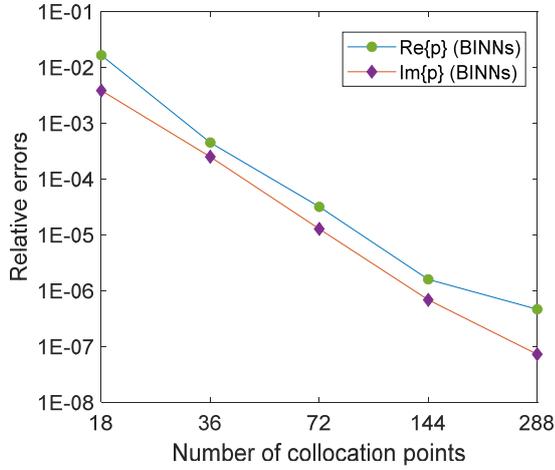 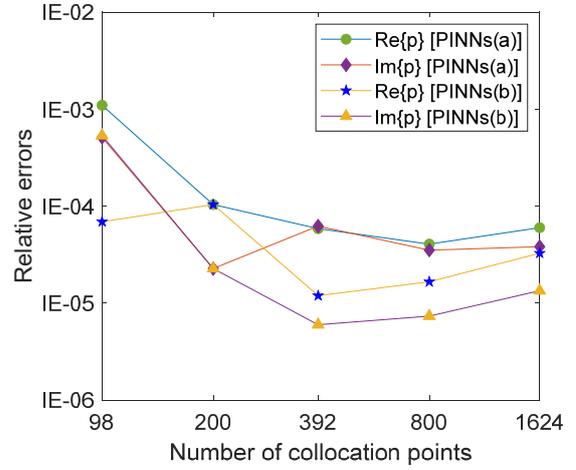

**Figure 5.** Convergence curves of pressures by the BINNs with different number of collocation points

**Figure 6.** Convergence curves of pressures by the PINNs with different number of collocation points

Case 2: Mixed BCs

The mixed BCs are taken into account in this particular case. As depicted in Fig. 3, the left, upper and lower boundaries of the domain are assumed to be rigid, while the right boundary is subjected to a specific condition as

$$p(3, x'_2) = \sin x'_2 + \mathrm{i} \cos x'_2, \ (3, x'_2) \in \Gamma \tag{23}$$

The analytical solution for the case is not available.

The wave number is assumed to be $k = 2$ m$^{-1}$. Both the BINNs and the PINNs are employed for the numerical simulation of this case to make a comparison. The activation functions of the neural networks remain the same as in case 1, and the training process for optimization stops after 10000 iterations. The BINNs uses 288 collocation points and a single hidden layer with 20 neurons, while the PINNs uses 1624 collocation points and three hidden layer, each with 25 neurons. Fig. 7 displays the numerical results of the pressures in the entire computational domain. As observed from the figure, the numerical results obtained by the BINNs show good agreement with those calculated by the PINNs.

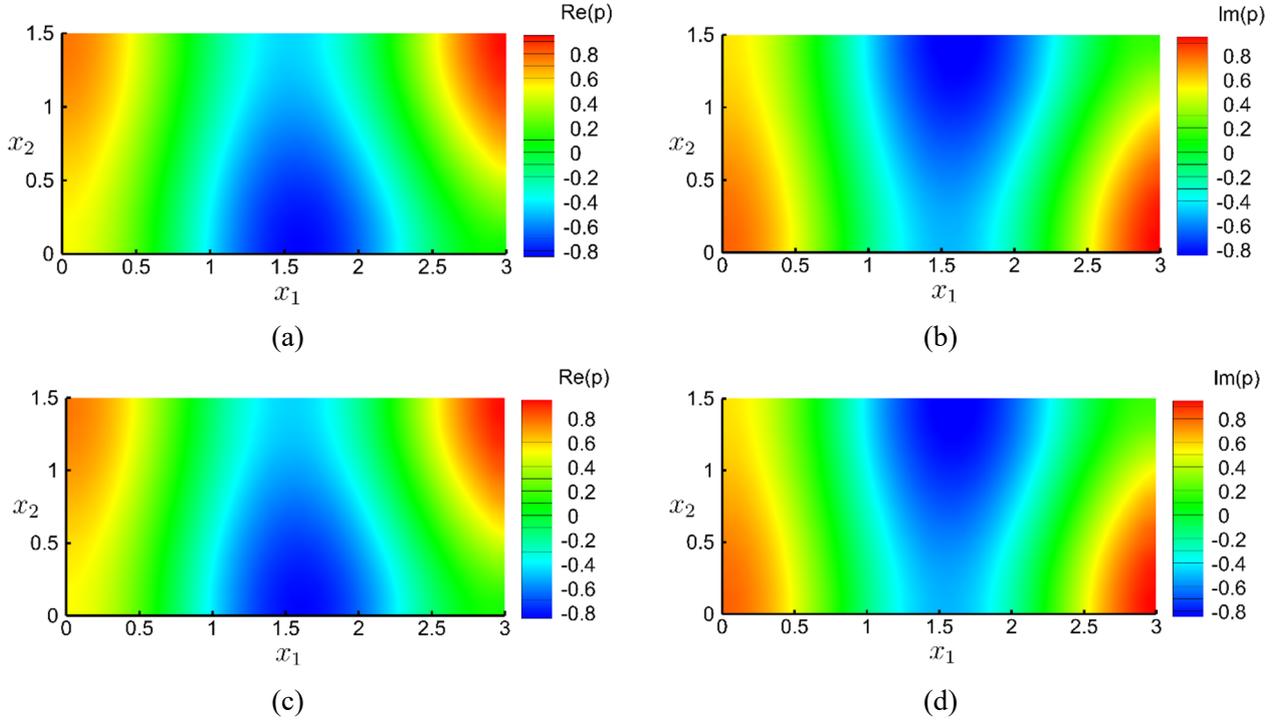

**Figure 7.** Numerical results of pressures in the rectangle domain: a) real component (the BINNs); b) imaginary component (the BINNs); c) real component (the PINNs); d) imaginary component (the PINNs)

*4.2. Acoustic radiation of an infinite pulsating cylinder*

The second example focuses on the analysis of acoustic radiation from an infinite pulsating cylinder. The cylinder has a radius of $R = 1\,\text{m}$ and its center is located at $(0,0)$. The boundary of the structure has a normal velocity amplitude of $\bar{v} = 1\,\text{m/s}$. The analytical solution for the pressure can be determined as

$$p(r) = i\rho c \bar{v} \frac{H_0^1(kr)}{H_1^1(kR)}, \quad r \geq R. \tag{24}$$

where $H_i^1 (i = 0,1)$ denote the $i$-th order Hankel function of the first kind. The medium for the propagation of acoustic waves is assumed to be air, with a density of $\rho = 1.2\,\text{kg/m}^3$ and a wave speed of $c = 341\,\text{m/s}$.

In this simulation, the wave number $k = 1\,\text{m}^{-1}$ is selected. The BINNs employs neural networks consisting of two hidden layers, each comprising 10 neurons. The training process for optimization terminates after 2000 iterations. The present approach utilizes 150 collocation points on the boundary. The chosen activation function is "Swish", as specified in Table 1. Calculated points are distributed within a domain $\{(x_1, x_2) | \sqrt{x_1^2 + x_2^2} > 1, -5 < x_1, x_2 < -5\}$. Fig. 8 presents the contour plots of relative errors for the real and imaginary components of pressures at the calculated points, as evaluated by the BINNs. It is evident that the present approach yields satisfactory numerical results.

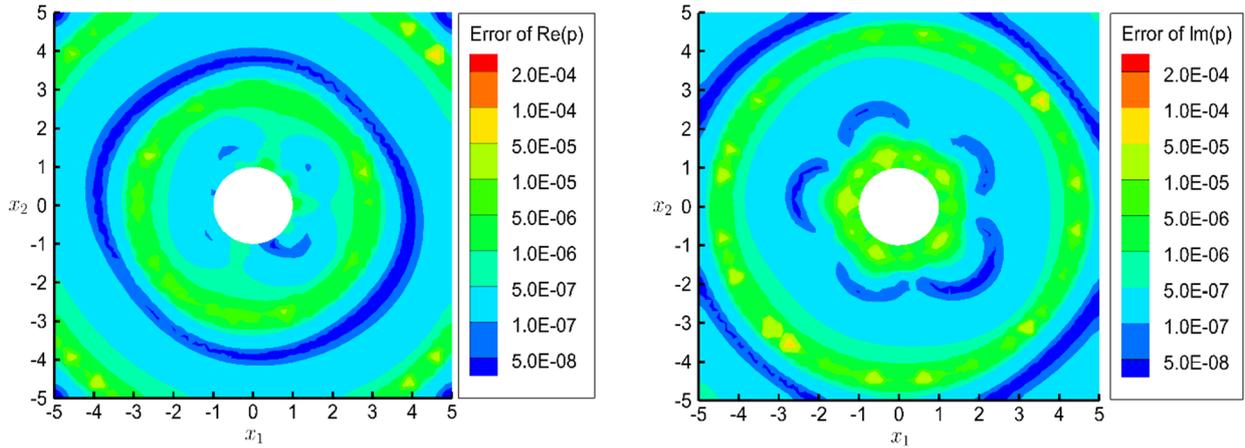

**Figure 8.** Numerical results of pressures calculated by the BINNs

Maintaining the previous settings unaltered, we proceed to validate the impact of various activation functions listed in Table 1 on the developed method. Table 3 shows the numerical errors of pressures in domain $\{(x_1, x_2) | \sqrt{x_1^2 + x_2^2} > 1, -5 < x_1, x_2 < -5\}$, along with the CPU time and the final values of *Loss*, obtained using the BINNs with different activation functions. From the table, it indicates that the choice of activation functions has minimal effect on the precision, convergence process of the loss function, and the efficiency of the BINNs.

Table 3. Impact of various activation functions on the BINNs

| Activation functions | Arctan | Sigmoid | Swish | Softplus | Tanh |
|---|---|---|---|---|---|
| Error of Re{$p$} | 1.10E-06 | 5.44E-06 | 9.97E-07 | 4.59E-07 | 1.02E-06 |
| Error of Im{$p$} | 1.20E-06 | 6.72E-06 | 1.02E-06 | 6.67E-07 | 2.19E-06 |
| Final value of $Loss$ | 7.15E-10 | 1.47E-09 | 4.87E-10 | 8.95E-11 | 2.95E-09 |
| CPU time (s) | 21.6 | 23.5 | 21.9 | 22.0 | 21.8 |

*4.3. Acoustic scattering of an infinite rigid cylinder*

As the last numerical example, we consider an acoustic scattering phenomenon. A plane incident wave, with an amplitude of unity, travels along the positive *x*-axis and impinges on an infinite rigid cylinder centered at point (0, 0) with a radius of $R = 1\,\text{m}$. The analytical solution of scattering field

$$p(r,\theta) = -\sum_{n=0}^{\infty} \varepsilon_n \mathrm{i}^n \frac{J'_n(kR)}{H_n^{1'}(kR)} H_n^1(kr) \cos(n\theta), \ r \geq R \tag{25}$$

where $J_n$ denotes the *n*-th order Bessel function, $H_n^1$ represents the *n*-th order Hankel function of the first kind, $\theta = 0$ along the positive *x*-axis, and $\varepsilon_n$ is the Neumann symbol expressed as

$$\varepsilon_n = \begin{cases} = 1, & n = 0, \\ = 2, & n \geq 1. \end{cases} \tag{26}$$

A neural network with a configuration of two hidden layers, each consisting of 20 neurons, is utilized for the numerical implementation of the BINNs. The wave number is set to $k = 0.5\,\text{m}^{-1}$, and a total of 90 collocation points are distributed on the boundary. The activation function is set to $\sigma(z) = z/(1+e^{-z})$. Two loss functions, specifically Eq. (19) and Eq. (20), are reconsidered and incorporated into the BINNs for analyzing acoustic fields in unbounded domains. Fig. 9 depicts the convergence behavior of two designated loss functions, namely $Loss$ and $Loss_{BC}$, as the iterations progress from 1 to 2000, with measurements taken every 50 iterations. Once again, it is demonstrated that $Loss$ has a better convergence performance when compared to $Loss_{BC}$.

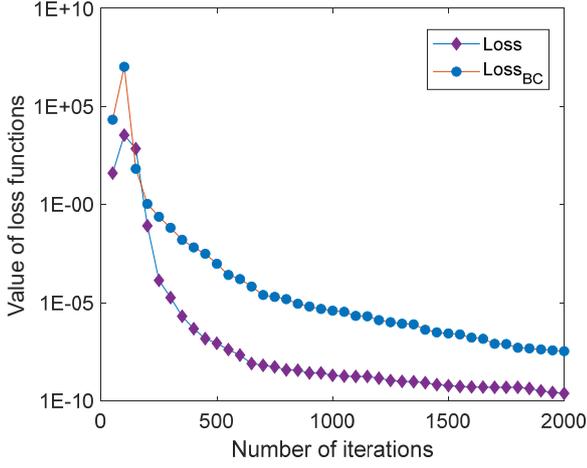 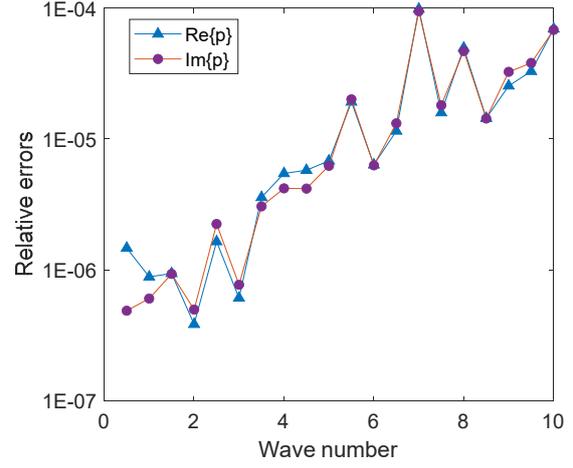

**Figure 9.** Convergence process of loss functions

**Figure 10.** Variations of relative errors of pressures with different wave numbers

The number of collocation points on the boundary is adjusted to 300, and the training process for optimization is conducted over 5000 iterations. All other settings remain unchanged from the previous configuration. Employing the BINNs with $Loss$, Fig. 10 displays the relative errors of pressures in domain $\{(x_1, x_2) | 1 < \sqrt{x_1^2 + x_2^2} < 2\}$ across varying wave numbers ranging from $0.5 \text{ m}^{-1}$ to $10 \text{ m}^{-1}$. As we can observe in Fig. 10, the developed method obtains the accurate numerical results for different wave numbers. Fig. 11 presents the relative errors of pressures at all calculated points within domain $\{(x_1, x_2) | 1 < \sqrt{x_1^2 + x_2^2} < 2\}$, considering a wave number of $k = 5 \text{ m}^{-1}$. It can be observed that maximum relative error of both the real and imaginary parts of pressures at these calculated points is below 5E-003.

These numerical results obtained using the BINNs further illustrate the competitiveness of the proposed method in simulating acoustic fields in unbounded domains, surpassing the traditional PINNs.

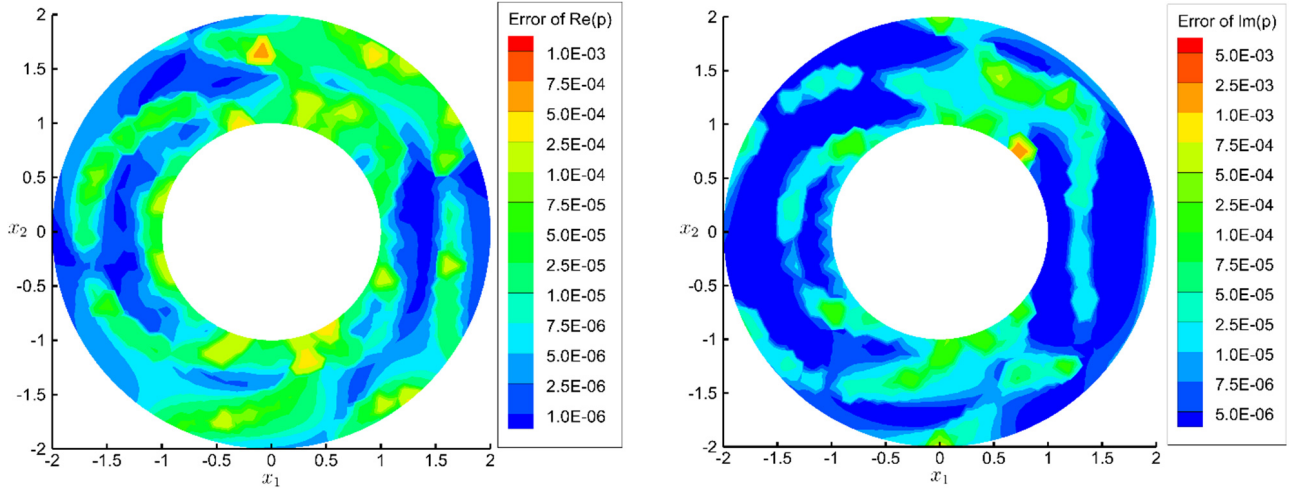

**Figure 11.** Numerical errors of pressures calculated by the BINNs for $k = 5 \text{ m}^{-1}$

## 5. Conclusions

The BINNs is proposed in this paper as a numerical approach for analyzing acoustic fields in both bounded and unbounded domains. Unlike the traditional PINNs that combines the governing equation with neural architectures, the proposed method integrates the BIEs and neural networks. Through numerical experiments on various benchmark examples, the BINNs exhibits high accuracy and rapid convergence. Several notable advantages of the BINNs over the traditional PINNs in the context of acoustic radiation and scattering can be summarized as follows:

1) The BINNs only require the coordinates of "boundary" collocation points as input data for the neural networks. The benefit of this is that the method is particularly well-suited for numerical simulations of problems in unbounded domains.

2) The loss function in the BINNs, as defined in Eq. (19), is not a composite form. Therefore, there is no need to consider special techniques to balance the influence between different terms, as described in Eq. (20) or the loss function used in the PINNs. The numerical results also demonstrate the fast convergence of the loss function.

3) To achieve comparable precision in pressure calculations, the BINNs requires significantly fewer collocation points and hidden layers/neurons compared to the PINNs. As a result, the BINNs exhibits higher computational efficiency.

4) The BINNs has higher precision attributed to the semi-analytic characteristic of the BIEs, as evident from the numerical errors of acoustic pressures obtained using this method.

The present approach is introduced to address relatively simple acoustic problems, and several conclusions are summarized. In the future, we aim to extend the application of BINNs to structural-acoustic sensitivity analysis.

## Acknowledgements

The work described in this paper was supported by the Natural Science Foundation of Shandong Province of China (Grant No. ZR2022YQ06), the Development Plan of Youth Innovation Team in Colleges and Universities of Shandong Province (Grant No. 2022KJ140), the National Natural Science Foundation of China (Grant No. 11802165), and the China Postdoctoral Science Foundation (Grant No. 2019M650158).

## References

[1] Simpson RN, Scott MA, Taus M, Thomas DC, Lian H. Acoustic isogeometric boundary element analysis. Computer Methods in Applied Mechanics and Engineering 2014;269: 265-90.
[2] Zheng C, Matsumoto T, Takahashi T, Chen H. A wideband fast multipole boundary element method for three dimensional acoustic shape sensitivity analysis based on direct differentiation method. Engineering Analysis with Boundary Elements 2012;36(3): 361-71.
[3] Chen L, Zheng C, Chen H. A wideband FMBEM for 2D acoustic design sensitivity analysis based on direct differentiation method. Computational Mechanics 2013;52: 631-48.
[4] Jordan MI, Mitchell TM. Machine learning: Trends, perspectives, and prospects. Science 2015;349(6245): 255-60.
[5] Ruthotto L, Haber E. Deep neural networks motivated by partial differential equations. Journal of Mathematical Imaging and Vision 2020;62: 352-64.
[6] Zeng S, Zhang Z, Zou Q. Adaptive deep neural networks methods for high-dimensional partial differential equations. Journal of Computational Physics 2022;463: 111232.

[7] Raissi M, Perdikaris P, Karniadakis GE. Physics-informed neural networks: A deep learning framework for solving forward and inverse problems involving nonlinear partial differential equations. Journal of Computational physics 2019;378: 686-707.

[8] Yang L, Meng X, Karniadakis GE. B-PINNs: Bayesian physics-informed neural networks for forward and inverse PDE problems with noisy data. Journal of Computational Physics 2021;425: 109913.

[9] Gu Y, Zhang C, Zhang P, Golub MV, Yu B. Enriched physics-informed neural networks for 2D in-plane crack analysis: Theory and MATLAB code. International Journal of Solids and Structures 2023;276: 112321.

[10] Sirignano J, Spiliopoulos K. DGM: A deep learning algorithm for solving partial differential equations. Journal of computational physics 2018;375: 1339-64.

[11] Saporito YF, Zhang Z. Path-dependent deep Galerkin method: A neural network approach to solve path-dependent partial differential equations. SIAM Journal on Financial Mathematics 2021;12(3): 912-40.

[12] Yu B. The deep Ritz method: a deep learning-based numerical algorithm for solving variational problems. Communications in Mathematics and Statistics 2018;6(1): 1-12.

[13] Gu Y, Ng MK. Deep Ritz Method for the Spectral Fractional Laplacian Equation Using the Caffarelli--Silvestre Extension. SIAM Journal on Scientific Computing 2022;44(4): A2018-A36.

[14] Song C, Alkhalifah T, Waheed UB. Solving the frequency-domain acoustic VTI wave equation using physics-informed neural networks. Geophysical Journal International 2021;225(2): 846-59.

[15] Zhang Y, Zhu X, Gao J. Seismic inversion based on acoustic wave equations using physics-informed neural network. IEEE Transactions on Geoscience and Remote Sensing 2023;61: 1-11.

[16] Lin G, Hu P, Chen F, Chen X, Chen J, Wang J, et al. Binet: learning to solve partial differential equations with boundary integral networks. arXiv preprint arXiv:211000352 2021.

[17] Zhang H, Anitescu C, Bordas S, Rabczuk T, Atroshchenko E. Artificial neural network methods for boundary integral equations. TechRxiv Preprint techrxiv20164769v1 2022.

[18] Qu W, Chen W, Zheng C. Diagonal form fast multipole singular boundary method applied to the solution of high-frequency acoustic radiation and scattering. International Journal for Numerical Methods in Engineering 2017;111(9): 803-15.

[19] Guiggiani M, Casalini P. Direct computation of Cauchy principal value integrals in advanced boundary elements. International Journal for Numerical Methods in Engineering 1987;24(9): 1711-20.

[20] Qu W, Chen W, Fu Z. Solutions of 2D and 3D non-homogeneous potential problems by using a boundary element-collocation method. Engineering Analysis with Boundary Elements 2015;60: 2-9.